\documentclass[12pt,leqno,draft]{article}
\usepackage{amsfonts}
\usepackage{amsmath, amsthm, amsfonts, amssymb, color}
\usepackage{mathrsfs}
\usepackage{color}
\usepackage[notcite,notref]{showkeys}

\newtheorem{thm}{Theorem}[section]

\newtheorem{exa}[thm]{Example}
\newtheorem{rem}[thm]{Remark}
\theoremstyle{definition}
\newtheorem{defn}{Definition}[section]
\newcommand{\scr}[1]{\mathscr #1}
\definecolor{wco}{rgb}{0.5,0.2,0.3}

\numberwithin{equation}{section} \theoremstyle{remark}

\newcommand{\ua}{\uparrow}

\title{{\bf Weak Solution and Invariant Probability Measure for McKean-Vlasov SDEs with Integrable Drifts}\footnote{Supported in
 part by  NNSFC (11801406).} }
\author{
{\bf   Xing Huang $^{a)}$,  Shen Wang $^{a)}$, Fen-Fen Yang $^{b)}$ }\\
\footnotesize{ a)Center for Applied Mathematics, Tianjin
University, Tianjin 300072, China}\\
\footnotesize{  xinghuang@tju.edu.cn, wswangshen@tju.edu.cn}\\
 \footnotesize{ b)Department of Mathematics, Shanghai University, Shanghai 200444, China}\\
\footnotesize{ yangfenfen@shu.edu.cn }}
\begin{document}
\allowdisplaybreaks
\def\R{\mathbb R}  \def\ff{\frac} \def\ss{\sqrt} \def\B{\mathbf
B} \def\W{\mathbb W}
\def\N{\mathbb N} \def\kk{\eta} \def\m{{\bf m}}
\def\ee{\varepsilon}\def\ddd{D^*}
\def\dd{\delta} \def\DD{\Delta} \def\vv{\varepsilon} \def\rr{\rho}
\def\<{\langle} \def\>{\rangle} \def\GG{\Gamma} \def\gg{\gamma}
  \def\nn{\nabla} \def\pp{\partial} \def\E{\mathbb E}
\def\d{\text{\rm{d}}} \def\bb{\beta} \def\aa{\alpha} \def\D{\scr D}
  \def\si{\sigma} \def\ess{\text{\rm{ess}}}
\def\beg{\begin} \def\beq{\begin{equation}}  \def\F{\scr F}
\def\Ric{\text{\rm{Ric}}} \def\Hess{\text{\rm{Hess}}}
\def\e{\text{\rm{e}}} \def\ua{\underline a} \def\OO{\Omega}  \def\oo{\omega}
 \def\tt{\tilde} \def\Ric{\text{\rm{Ric}}}
\def\cut{\text{\rm{cut}}} \def\P{\mathbb P} \def\ifn{I_n(f^{\bigotimes n})}
\def\C{\scr C}      \def\aaa{\mathbf{r}}     \def\r{r}
\def\gap{\text{\rm{gap}}} \def\prr{\pi_{{\bf m},\varrho}}  \def\r{\mathbf r}
\def\Z{\mathbb Z} \def\vrr{\varrho}
\def\L{\scr L}\def\Tt{\tt} \def\TT{\tt}\def\II{\mathbb I}
\def\i{{\rm in}}\def\Sect{{\rm Sect}}  \def\H{\mathbb H}
\def\M{\scr M}\def\Q{\mathbb Q} \def\texto{\text{o}}
\def\Rank{{\rm Rank}} \def\B{\scr B} \def\i{{\rm i}} \def\HR{\hat{\R}^d}
\def\to{\rightarrow}\def\l{\ell}\def\iint{\int}
\def\EE{\scr E}\def\Cut{{\rm Cut}}
\def\A{\scr A} \def\Lip{{\rm Lip}}
\def\BB{\scr B}\def\Ent{{\rm Ent}}\def\L{\scr L}
\def\R{\mathbb R}  \def\ff{\frac} \def\ss{\sqrt} \def\B{\mathbf
B}
\def\N{\mathbb N} \def\kk{\eta} \def\m{{\bf m}}
\def\dd{\delta} \def\DD{\Delta} \def\vv{\varepsilon} \def\rr{\rho}
\def\<{\langle} \def\>{\rangle} \def\GG{\Gamma} \def\gg{\gamma}
  \def\nn{\nabla} \def\pp{\partial} \def\E{\mathbb E}
\def\d{\text{\rm{d}}} \def\bb{\beta} \def\aa{\alpha} \def\D{\scr D}
  \def\si{\sigma} \def\ess{\text{\rm{ess}}}
\def\beg{\begin} \def\beq{\begin{equation}}  \def\F{\scr F}
\def\Ric{\text{\rm{Ric}}} \def\Hess{\text{\rm{Hess}}}
\def\e{\text{\rm{e}}} \def\ua{\underline a} \def\OO{\Omega}  \def\oo{\omega}
 \def\tt{\tilde} \def\Ric{\text{\rm{Ric}}}
\def\cut{\text{\rm{cut}}} \def\P{\mathbb P} \def\ifn{I_n(f^{\bigotimes n})}
\def\C{\scr C}      \def\aaa{\mathbf{r}}     \def\r{r}
\def\gap{\text{\rm{gap}}} \def\prr{\pi_{{\bf m},\varrho}}  \def\r{\mathbf r}
\def\Z{\mathbb Z} \def\vrr{\varrho}
\def\L{\scr L}\def\Tt{\tt} \def\TT{\tt}\def\II{\mathbb I}
\def\i{{\rm in}}\def\Sect{{\rm Sect}}  \def\H{\mathbb H}
\def\M{\scr M}\def\Q{\mathbb Q} \def\texto{\text{o}}
\def\Rank{{\rm Rank}} \def\B{\scr B} \def\i{{\rm i}} \def\HR{\hat{\R}^d}
\def\to{\rightarrow}\def\l{\ell}
\def\8{\infty}\def\I{1}\def\U{\scr U}
\maketitle

\begin{abstract}
In this paper, by utilizing Wang's Harnack inequality with power and the Banach fixed point theorem, the weak well-posedness for
McKean-Vlasov   SDEs with integrable drift is investigated. In addition,
using the decoupled method, some regularity such as relative entropy and Sobolev's estimate of invariant probability measure are proved. Finally, by Banach's fixed theorem, the existence and uniqueness of invariant probability measure for symmetric McKean-Vlasov SDEs and stochastic Hamiltonian system with integrable drifts are obtained.
\end{abstract} \noindent
 AMS subject Classification:\  60H10, 60G44.   \\
\noindent
 Keywords: McKean-Vlasov SDEs, Total variation distance, Weak solution, Invariant probability measure, Integrable condition.
 \vskip 2cm

\section{Introduction}
Invariant probability measure is the equilibrium state in physics. There are plentiful results on the invariant probability measure for linear semigroup $P_t$ associated to classical diffusion process in $\R^d$:
 $$\d X_t=b(X_t)\d t+\sigma(X_t)\d W_t,$$
 the infinitesimal generator of which is defined as $$L=\frac{1}{2}\mathrm{Tr}(\sigma\sigma^\ast\nabla ^2)+b\cdot\nabla.$$
The existence of invariant probability measure can be studied by investigating the tightness of the sequence of probability measures
$$\frac{1}{n}\int_0^n P^\ast_t\delta_x\d t,\ \ n\geq 1,$$
see \cite{DZ1}.
Meanwhile, a useful sufficient condition to obtain the existence of invariant probability measure is Lyapunov's condition, i.e.
\begin{align*}LW_1\leq C-W_2
 \end{align*}
for some positive function $W_1\in C^2(\R^d)$, positive compact function $W_2$ and some constant $C>0$ can derive that $P_t$ has an invariant probability measure $\mu$ satisfying $\mu(W_2)\leq C$, see \cite{BKR,BRW,BR}.

For the uniqueness, the classical principle is strong Feller property together with irreducibility, see \cite[Theorem 4.2.1]{DZ1}. Moreover, by Wang's Harnack inequality \cite[Theorem 1.4.1]{Wbook}, the uniqueness can also
be ensured.
Furthermore, using couplings or generalized couplings, \cite{S} proved the uniqueness of the invariant measures. One can also refer to \cite{BRS1,BRS2,BRS3,B} for conditions on the existence and uniqueness of invariant probability measure by Lyapunov function $V\in C^2(\R^d)$ with $\lim_{|x|\to\infty}V(x)=\infty$ and constants $C,R>0$ such that
\begin{align*}
L V\leq -C,\ \ |x|\geq R.
\end{align*}
Recently, in \cite{Wp}, the existence and uniqueness as well as the regularity such as relative entropy and Sobolev's estimate are derived by hyperboundedness or log-Sobolev's inequality.

However, all the above methods are invalid  to obtain the existence and uniqueness for distribution dependent SDEs(McKean-Vlasov SDEs or mean field SDEs), where the associated semigroup $P_t^\ast$ is nonlinear. In \cite{FYW1}, Wang obtained the existence and uniqueness of invariant probability measure and the exponential ergodicity in Wasserstein distance by the method of basic coupling(\cite[Definition 2.4]{CL}), see \cite{HRW} for the path-distribution dependent case and \cite{LMW} for the McKean-Vlasov SDEs with L\'{e}vy noise. Quite recently, \cite{Zhang} investigated the existence of invariant probability measure by Schauder's fixed point theorem, see also \cite{BSY} for the existence of invariant probability measure of functional McKean-Vlasov SDEs by Kakutani's fixed point theorem. In addition, \cite{W21} proved the existence and uniqueness of invariant probability measure for (reflecting) McKean-Vlasov SDEs by exponential ergodicity of the decoupled SDEs and Banach's fixed point theorem. For more results, one can see \cite{D,DT,FZ} and references therein. What's more, by using log-Sobolev's inequality or Poinc\'{a}re's inequality for the invariant probability measure of decoupled SDE and the Banach fixed point theorem, \cite{BK,BR} investigated the existence and uniqueness of the solution to stationary nonlinear and non-degenerate Fokker-Planck-Kolmogorov equations.

In this paper, we will prove the weak well-posedness for McKean-Vlasov SDEs with drift being integrable in the spacial component with respect to some reference probability measure by Banach's fixed point theorem. The result extend the one in \cite{Wp}. Moreover, the regularity of invariant probability measure with integrable drift is also obtained by decoupled technique and the existed results in \cite{Wp}. Since the invariant probability measure of decoupled SDE is in general unknown when the drift is only assumed to be integrable, the conditions in  \cite{BK,BR} such as log-Sobolev's inequality and Poincar\'{e}'s inequality are not explicit. As a result, we only investigate the existence and uniqueness of the invariant probability measure for symmetric and non-degenerate McKean-Vlasov SDEs as well as distribution dependent stochastic Hamiltonian system, where the drift is assumed to be of gradient form and  integrable in the spacial component with respect to some reference probability measure.

Let $\scr P(\R^d)$ be the space of all probability measures on $\R^d$ equipped with the weak topology.
 Consider the following   distribution dependent SDE on $\R^d$:
\beq\label{E1} \d X_t= \{Z_0(X_t)+\sigma(X_t)Z(X_t, \L_{X_t})\}\d t +\si(X_t)\d W_t,\end{equation}
 where $(W_t)_{t\geq 0}$ is an $n$-dimensional Brownian motion on a complete filtration probability space $(\OO,\F,\{\F_t\}_{t\geq 0},\P)$, $\L_{X_t}$ is the law of $X_t$,
$$Z: \R^d\times \scr P(\R^d)\to \R^n,\ \ Z_0: \R^d\to \R^d,\ \  \si: \R^d\to \R^d\otimes\R^n$$ are measurable. Compared with \cite{Wp}, $Z$ can depend on the distribution of the solution, see $(A_b)$ below for the condition of $Z$ on the measure component.
When a different probability measure $\tt\P$ is concerned, we use $\L_\xi|\tt \P$ to denote the law of a random variable $\xi$ under the probability $\tt\P$, and use $\E_{\tt\P}$ to stand for  the expectation under $\tt\P$.
\begin{defn}
\begin{enumerate}
\item[(1)] An adapted continuous process $(X_t)_{t\geq 0}$ on $\R^d$ is called a solution of \eqref{E1}, if $X_0$ is $\F_0$-measurable,
\beq\label{A1} \E\int_0^T\big\{|Z_0(X_t)|+|\sigma(X_t) Z(X_t,\L_{X_t})|+\|\si(X_t)\|^2\big\}\d t<\infty,\ \ T>0,\end{equation}
and $\P$-a.s.
\beq\label{A2} \begin{split}X_t&= X_0+\int_0^{t}Z_0(X_s)\d s+ \int_0^{t} \sigma(X_s)Z(X_s,\L_{X_s})\d s +\int_0^t \si(X_s)\d W_s,\ \ t\geq 0.
\end{split}\end{equation}

\item[(2)]For any $\mu_0\in \scr P(\R^d)$,   $((\tt X_t)_{t\geq 0},(\tt W_t)_{t\geq 0})$  is called a weak solution to \eqref{E1} starting at $\mu_0$, if
$(\tt W_t)_{t\geq 0}$ is an $n$-dimensional Brownian motion under a   complete filtration probability space $(\tt\OO,\tilde{\F},\{\tt\F_t\}_{t\geq 0},\tt\P),$ $(\tilde{X}_t)_{t\geq 0}$ is   a continuous $\tt\F_t$-adapted process on $\R^d$ with $\L_{\tt X_0}|\tt \P=\mu_0$ and $\tilde{X}_0\in \tilde{\F}_0$, and
\eqref{A1}-\eqref{A2} hold for $(\tt X, \tt W, \tt\P,\E_{\tt \P})$
replacing
 $(X, W, \P,\E).$
\item[(3)] We call \eqref{E1} weakly well-posed for an initial distribution $\mu_0$, if it has a weak solution starting at $\mu_0$ and any weak solution with the same initial distribution is equal in law.
\end{enumerate}
\end{defn}
For the well-posedness of distribution dependent SDEs with singular drifts, one can refer to \cite{BB,BBP,CA,CF,CR,HW,HW21,L,MV,RZ,ZG} and references within.

The remaining part of the paper is organized as follows:
In Section 2, we investigate the weak well-posedness of \eqref{E1} under integrable condition.
In Section 3, the regularity of invariant probability measure for McKean-Vlasov SDEs with integrable drift is presented. The existence and uniqueness of invariant probability measure are provided in Section 4.
  \section{Weak Well-posedness}
 For any $\mu,\nu\in\scr P(\R^d)$, the total variation distance between $\mu$ and $\nu$ is defined as
$$\|\mu-\nu\|_{TV}=2\sup_{A\in\scr B(\R^d)}|\mu(A)-\nu(A)|=\sup_{\|f\|_\infty\leq 1}|\mu(f)-\nu(f)|.$$
To obtain the weak well-posedness of \eqref{E1}, we make the following assumptions, see \cite{Wp} for more details.
\begin{enumerate}
\item[$({\bf A})$] The reference SDE
\beq\label{EM} \d X_t= Z_0(X_t)\d t+\si(X_t)\d W_t\end{equation}
is strongly well-posed and has a unique invariant probability measure $\mu^0$.
\end{enumerate}
\begin{enumerate}
\item[$(A_Z)$] 
    There exist  constants $\varepsilon>0, K_Z>0$ such that
\begin{align}\label{ne}\mu^0(\e^{\varepsilon|Z(\cdot,\mu^0)|^2})<\infty,\end{align}
and
 \begin{equation}\label{cono}\beg{split}
&|Z(x,\gamma)-Z(x,\tilde{\gamma})|\le K_Z\|\gamma-\tilde{\gamma}\|_{TV},\ \ x\in \R^d, \gamma,\tilde{\gamma}\in \scr P(\R^d).
\end{split}
\end{equation}
\end{enumerate}
Let $\pi_t$ be the projection map from $C([0,\infty);\R^d)$ to $\R^d$, i.e. $$\pi_t(w)=w_t, \ \ w\in C([0,\infty);\R^d).$$
For any $\gamma\in \scr P(\R^d)$, we will prove that \eqref{E1} has a unique weak solution with initial distribution $\gamma$ and the distribution of the solution $\P^\gamma$ satisfying
\begin{align}\label{Pg}\P^\gamma\left(w\in C([0,\infty);\R^d), \int_0^t |Z(w_s,\P^\gamma\circ \pi_s^{-1})|^2\d s<\infty,\  \ t\geq0\right)=1.
\end{align}
Firstly, repeating the proof of \cite[Theorem 2.1]{Wp}, we can easily extend {it} to the time dependent case below.
Consider
\beq\label{E1'} \d X_t=\{Z_0(X_t)+\sigma(X_t)\tilde{Z}_t(X_t)\}\d t +\si(X_t)\d W_t,\end{equation}
here $\tilde{Z}:[0,\infty)\times \R^d\to \R^n$ is measurable.
For any $\gamma\in\scr P(\R^d)$, let $X_t^\gamma$ be the solution to \eqref{EM} with initial distribution $\gamma$.
\begin{thm}\label{zi} Assume $({\bf A})$ and that there exists a constant $\varepsilon>0$ such that $$\|\e^{\varepsilon|\tilde{Z}|^2} \|_{L^{\infty}([0,t];L^1(\mu^0))}<\infty, \ \ t>0.$$
\begin{enumerate}
\item[(i)] For any $\gamma\in\scr P(\R^d)$ with $\gamma\ll\mu^0$,
\eqref{E1'} has a unique weak solution starting at $\gamma$ and the distribution $\P^\gamma$ satisfying
\begin{align}\label{Pgo}\P^\gamma\left(w\in C([0,\infty);\R^d), \int_0^t |\tilde{Z}_s(w_s)|^2\d s<\infty,\  \ t\geq0\right)=1.
\end{align}
Moreover, it holds
$$P_t^{\tilde{Z}}f(\gamma):=[\P^\gamma\circ\pi_t^{-1}](f)=\E (R^\gamma(t)f(X^\gamma_t)),\ \ t\geq 0, f\in\scr B_b(\R^d)$$
where $R^\gamma(t)$ is a martingale given by
$$R^\gamma(t)=\exp\left\{\int_0^t\<\tilde{Z}_s(X_s^\gamma),\d W_s\>-\frac{1}{2}\int_0^t|\tilde{Z}_s(X_s^\gamma)|^2\d s\right\},\ \ t\geq 0.$$
\item[(ii)]
If in addition, there exists $p>1$ such that
\begin{align}\label{HI}
(P_t^0|f|)^p(z)\leq (P_t^0|f|^p)(\bar{z})\e^{\Phi_p(t,z,\bar{z})},\ \ f\in \B_b(\R^d), z,\bar{z}\in\R^d,t>0
\end{align}  with
\begin{align}\label{E0p}
\int_0^t\left\{\mu^0(\e^{-\Phi_p(s,z,\cdot)})\right\}^{-\frac{1}{p}}\d s<\infty, \ \ t>0, z\in\R^d,
\end{align}
then the assertion in (i) holds for any $\gamma\in\scr P(\R^d)$.
\end{enumerate}
\end{thm}
\begin{proof} Since the proof can be completely the same with that of \cite[Theorem 2.1]{Wp}, we omit it here.
\end{proof}
The main result in this section is the following theorem.
\begin{thm}\label{T1.1} Assume $({\bf A})$ and $(A_Z)$.
\begin{enumerate}
\item[(i)] For any $\gamma\in\scr P(\R^d)$ with $\gamma\ll\mu^0$,
\eqref{E1'} has a unique weak solution starting at $\gamma$ and satisfying \eqref{Pg}.
\item[(ii)] If in addition, $P_t^0$ satisfies \eqref{HI} and \eqref{E0p},
then for any $\gamma\in\scr P(\R^d)$, then (i) holds for any $\gamma\in\scr \P(\R^d)$.
\end{enumerate}
\end{thm}
\begin{rem}\label{ING} Compared with the localized integrable condition of the drift on the spacial component for the weak well-posedness in \cite[Theorem 3.9]{ZG}, the drift $Z$ in Theorem \ref{T1.1} is allowed to be of some growth. For instance, taking $n=d$, $Z_0=-x$, $\sigma=\sqrt{2}I_{d\times d}$, if $Z$ satisfies
$$|Z(x,\mu_0)|\leq c(1+|x|),$$
then \eqref{ne} holds with $c^2\varepsilon<\frac{1}{2}$.
Moreover, \eqref{E1} can be degenerate, see \cite[Example 4.3]{Wp} and Example \ref{SHS} below.
\end{rem}
To prove Theorem \ref{T1.1}, it is sufficient to prove that for any $T>0$, \eqref{E1} is weakly well-posed on $[0,T]$. So, we fix $T>0$ in the following.
For any
$\gg\in \scr P(\R^d)$, $\mu\in \scr B([0,T];\scr P(\R^d))$, consider
  \beq\label{E2} \d X_t=\{Z_0(X_t)+\sigma(X_t)Z(X_t, \mu_t)\}\d t +\si(X_t)\d W_t\end{equation}
  with initial distribution $\gamma$.
\begin{proof}[Proof of Theorem \ref{T1.1}] Since the proofs of (i) and (ii) are completely the same, we only prove (ii).

Note that  ({\bf$A_Z$}) implies that for any $t\in[0,T]$,
\begin{align}\label{INE}
\|\e^{\frac{\varepsilon}{2}|Z(\cdot,\mu_\cdot)|^2}\|_{L^\infty([0,t];L^1(\mu^0))}\leq \mu^0(\e^{\varepsilon|Z(\cdot,\mu^0)|^2+4\varepsilon K_Z^2})=\mu^0(\e^{\varepsilon|Z(\cdot,\mu^0)|^2})\e^{4\varepsilon K_Z^2}<\infty.
\end{align}
This together with $({\bf A})$, \eqref{HI}, \eqref{E0p}
and Theorem \ref{zi}(ii) for $\tilde{Z}^\mu_t=Z(\cdot,\mu_t)$ yields that
for any $\gamma\in\scr P(\R^d)$,
\eqref{E2}  has a unique weak solution 
starting at $\gamma$ and the distribution $\P^{\gamma,\mu}$ satisfying
\begin{align}\label{PGo}\P^{\gamma,\mu}\left(w\in C([0,T];\R^d), \int_0^T |Z(w_s,\mu_s)|^2\d s<\infty\right)=1.
\end{align}
Define $\Phi_t^\gg(\mu)=\P^{\gamma,\mu}\circ\pi_t^{-1},\ \ t\in[0,T]$. Let
$$R^{\gamma,\mu}(t)=\exp\left\{\int_0^t\<Z(X_s^\gamma,\mu_s),\d W_s\>-\frac{1}{2}\int_0^t|Z(X_s^\gamma,\mu_s)|^2\d s\right\},\ \ t\in[0,T].$$
Theorem \ref{zi}(ii) also implies
$$\Phi_t^\gg(\mu)(f)=\E( R^{\gamma,\mu}(t)f(X_t^\gamma)),\ \ t\in[0,T].$$
 For $\nu\in \scr B([0,T];\scr P(\R^d))$, let $R^{\gamma,\nu}(t)$ be defined in the same way as $R^{\gamma,\mu}(t)$ with $\nu$ in place of $\mu$. Then it is easy to see that
 $$\Phi_t^\gg(\nu)(f)=\E( R^{\gamma,\nu}(t)f(X_t^\gamma))=\E\left( R^{\gamma,\mu}(t) \frac{R^{\gamma,\nu}(t)}{ R^{\gamma,\mu}(t)}f(X_t^\gamma)\right),\ \ t\in[0,T].$$
Let
$$W_t^\nu=W_t-\int_{0}^{t}Z(X_s^\gamma,\nu_s)\d s,\ \ \xi_t=  Z(X^\gamma_t, \nu_t)-Z(X^\gamma_t, \mu_t),\ \ t\in [0,T].$$
Since $R^{\gamma,\nu}(t)$ is a martingale due to Theorem \ref{zi} and \eqref{INE}, we conclude that $\{W_t^\nu\}_{t\in[0,T]}$ is an $n$-dimensional Brownian motion under probability measure $\Q=R^{\gamma,\nu}_T\P$.
So, applying Pinsker's inequality   \cite{CK,Pin}, we obtain
 \begin{align*}
\|\Phi_t^\gg(\nu)-\Phi_t^\gg(\mu)\|_{TV}^2
&\leq 2\E\left(R^{\gamma,\nu}_t\log \left(\frac{R^{\gamma,\nu}_t}{R^{\gamma,\mu}_t}\right)\right)\\
 &=  2\E \left\{ R^{\gamma,\nu}_t\left\{\int_0^t \<\xi_s,\d W_s^\nu\>+\frac{1}{2}\int_0^t |\xi_s|^2\d s\right\}\right\} \\
 &= \E \left[R^{\gamma,\nu}_t\int_0^t |Z(X^\mu_s, \mu_s)-Z(X^\mu_s, \nu_s)|^2\d s\right].
 \end{align*}
 This together with \eqref{cono} implies
\begin{equation}\label{ph}
\|\Phi_t^\gg(\nu)-\Phi_t^\gg(\mu)\|_{TV}^2\leq \int_0^tK_Z^2\|\mu_s-\nu_s\|^2_{TV}\d s.
\end{equation}
Take $\lambda= K_Z^2 $  and consider the   space
$E_{T}:= \{\mu\in\B([0,T];\scr P(\R^d)):\mu_0=\gg\}$ equipped with the complete metric
$$\rr(\nu,\mu):= \sup_{t\in [0,T]} \e^{-\lambda t}\|\nu_t-\mu_t\|_{TV}.$$
It follows from \eqref{ph} that
\beq\label{ESTM'}\beg{split}
\sup_{t\in[0,T]}\e^{-2\lambda t}\|\Phi_t^\gg(\nu)-\Phi_t^\gg(\mu)\|_{TV}^2
&\leq \sup_{t\in[0,T]}\int_0^tK_Z^2\e^{-2\lambda (t-s)}\e^{-2\lambda s}\|\mu_s-\nu_s\|^2_{TV}\d s\\
&\leq \sup_{s\in[0,T]}\e^{-2\lambda s}\|\mu_s-\nu_s\|^2_{TV}\sup_{t\in[0,T]}\int_0^tK_Z^2\e^{-2\lambda (t-s)}\d s\\
&\leq \frac{1}{2}\sup_{s\in[0,T]}\e^{-2\lambda s}\|\mu_s-\nu_s\|^2_{TV}.\end{split}
 \end{equation}
Then $\Phi^\gg$ is a strictly contractive map on $E_{T}$, so that the equation
\beq\label{PHIE} \Phi_t^\gg(\mu)= \mu_t,\ \ t\in [0,T]\end{equation}  has a unique solution $\mu\in E_{T}.$ This combined with \eqref{PGo} completes the proof.
\end{proof}
\section{Regularity of Invariant Probability Measure}
In this section, we consider the regularity of invariant probability measure of \eqref{E1} and a general result for the regularity will be presented.
Throughout this section, we assume that $\mu^Z=\rho \mu^0$ is an invariant probability measure of \eqref{E1}. It is clear that $\mu^Z$ is also an invariant probability measure of the following decoupled SDE:
\beq\label{EMI} \d X_t=\{ Z_0(X_t)+\sigma(X_t)Z(X_t, \mu^Z)\}\d t +\si(X_t)\d W_t.\end{equation}
Denote
$$\scr P_Z=\{\nu=\rho_\nu\mu^0:\nu \ \ \text{is an invariant probability measure of \eqref{EMI}}\}.$$
\begin{thm}\label{rp} Assume $({\bf A})$. If there exists $t>0$ such that $P_t^0$ has a strictly positive density with respect to $\mu^0$ and
\begin{align}\label{HYP}\|P_{t_0}^0\|_{L^2(\mu^0)\to L^{2p_0}(\mu^0)}<\infty\end{align}
for some $t_0>0$ and $p_0>1$. Then $(A_Z)$ with some $\varepsilon>\kappa_0:=\frac{t_0(3 p_0-1)}{2(p_0-1)}$ implies
$$\mu^0(\rho \log \rho)\leq \inf_{\tilde{\varepsilon}\in(\kappa_0,\varepsilon)}\frac{ {t_0(3 p_0-1)} \log \mu^0(\e^{ \varepsilon|Z(\cdot,\mu^0)|^2+4 K_Z^2\frac{\varepsilon\tilde{\varepsilon}} {\varepsilon-\tilde{\varepsilon}}})+4\tilde{\varepsilon} p_0\log\|P_{t_0}^0\|_{L^2(\mu^0)\to L^{2p_0}(\mu^0)}}{2\tilde{\varepsilon}(p_0-1)-t_0(3 p_0-1) }.$$
\end{thm}
\begin{proof}
 For any $\tilde{\varepsilon}\in(0,\varepsilon)$, it follows from \eqref{cono} that
\begin{align}\label{em}\nonumber\mu^0(\e^{\tilde{\varepsilon} |Z(\cdot,\mu^Z)|^2})&\leq  \mu^0(\e^{\tilde{\varepsilon} (|Z(\cdot,\mu^0)|+K_Z\|\mu^Z-\mu^0\|_{TV})^2})\\
&= \mu^0(\e^{ \tilde{\varepsilon}|Z(\cdot,\mu^0)|^2+2\tilde{\varepsilon}\frac{K_Z\|\mu^Z-\mu^0\|_{TV}} {\sqrt{\varepsilon-\tilde{\varepsilon}}}\sqrt{\varepsilon-\tilde{\varepsilon}}|Z(\cdot,\mu^0)|+\tilde{\varepsilon}K_Z^2\|\mu^Z-\mu^0\|_{TV}^2})\\
\nonumber&\leq \mu^0(\e^{ \varepsilon|Z(\cdot,\mu^0)|^2})\e^{ K_Z^2\frac{\varepsilon\tilde{\varepsilon}} {\varepsilon-\tilde{\varepsilon}}\|\mu^Z-\mu^0\|^2_{TV}}<\infty.
\end{align}
Since $\varepsilon >\kappa_0=\frac{t_0(3 p_0-1)}{2(p_0-1)}$, using \eqref{em} for $\tilde{\varepsilon}\in(\kappa_0,\varepsilon)$ and \cite[Theorem 3.1, Theorem 4.1]{Wp}, we conclude that $\mu^Z$ is the unique invariant probability measure of \eqref{EMI} in $\scr P_Z$. Furthermore, \cite[Theorem 4.1]{Wp} holds with $\varepsilon$ replaced by $\tilde{\varepsilon}\in(\kappa_0,\varepsilon)$, i.e.
\begin{align}\label{mm}\mu^0(\rho\log\rho)&\leq \frac{ {t_0(3 p_0-1)} \log \mu^0(\e^{\tilde{\varepsilon} |Z(\cdot,\mu^Z)|^2})+4\tilde{\varepsilon} p_0\log\|P_{t_0}^0\|_{L^2(\mu^0)\to L^{2p_0}(\mu^0)}}{2\tilde{\varepsilon}(p_0-1)-t_0(3 p_0-1)}.
\end{align}
Substituting \eqref{em} into \eqref{mm}, we complete the proof.
\end{proof}
To obtain the Sobolev estimate for $\rho$ by log-Sobolev's inequality of the reference SDE \eqref{EM}, let
$$Z_0=\frac{1}{2}\sum_{i,j=1}^d\{\partial_j(\sigma\sigma^\ast)_{ij}-(\sigma\sigma^\ast)_{ij}\partial _j V\}e_i$$
    for some $V\in C^2(\R^d)$.
Define $$\scr E_0(f,g)=\mu^0(\<\sigma^\ast \nabla f,\sigma^\ast \nabla g\>),\ \ f, g\in C_0^\infty(\R^d).$$ Let $H_\sigma^{1,2}(\mu^0)$ be the completion of $C_0^\infty(\R^d)$ under the norm
$$\sqrt{\scr E_1(f,f)}:=\{\mu^0(|f|^2+|\sigma^\ast\nabla f|^2)\}^{\frac{1}{2}}.$$
Then $(\scr E_0, H_\sigma^{1,2}(\mu^0))$ is a symmetric Dirichlet form on $L^2(\mu^0)$.

Moreover, we shall introduce the condition {\bf (H)} in \cite{Wp}:
\begin{enumerate}
\item[(\bf H)] Assume that $\mu^0(\d x)=\e^{-V}\d x$ is a probability measure. There exists $k\geq 2$ such that $\sigma\in C^k(\R^d,\R^d\otimes\R^n)$ and vector fields
    $$U_i=\sum_{j=1}^d\sigma_{ji}\partial _j,\ \ i=1,\cdots,n$$
    satisfy the H\"{o}rmander condition up to the $k$-th order of Lie brackets. Moreover, $1\in H_\sigma^{1,2}(\mu_0)$ with $\scr E_0(1,1)=0$, and defective log-Sobolev inequality
\begin{align}\label{dl}\mu^0(f^2\log f^2)\leq \kappa \mu^0(|\sigma^\ast\nabla f|^2)+\beta,\ \ f\in C_0^\infty(\R^d),\mu^0(f^2)=1
\end{align}
holds for some $\kappa>0$ and $\beta\geq 0$.
\end{enumerate}
One can refer to \cite{A,BE,BD,CW,DS,W} for more results on the log-Sobolev inequality.
\begin{thm}\label{SE} Assume {\bf (H)} and $(A_Z)$ for  some $\varepsilon >\kappa$. Then $\rho$ has a strictly positive continuous version satisfying $\log \rho$, $\rho ^{\frac{p}{2}}\in H^{1,2}_{\sigma}(\mu^0)$ for $p\in(1,\frac{\sqrt{\varepsilon}}{\sqrt{\kappa}})$ and
\begin{align*}
&\mu^0(|\sigma^\ast\nabla\log \rho|^2)\leq 4\mu^0((|Z(\cdot,\mu^0)|+2K_Z)^2)\\
&\mu^0(|\sigma^\ast\nabla\sqrt{\rho}|^2)\leq \inf_{\tilde{\varepsilon}\in(\kappa,\varepsilon)}\frac{1}{\tilde{\varepsilon}-\kappa}(\log \mu^0(\e^{ \varepsilon|Z(\cdot,\mu^0)|^2+4 K_Z^2\frac{\varepsilon\tilde{\varepsilon}} {\varepsilon-\tilde{\varepsilon}}})+\beta),\\
&\mu^0(|\sigma^\ast\nabla\rho^{\frac{p}{2}}|^2+\rho^p)\leq \inf_{\tilde{\varepsilon}\in(p^2\kappa,\varepsilon)}C_{p,\tilde{\varepsilon}}(\mu^0(\e^{ \varepsilon|Z(\cdot,\mu^0)|^2+4 K_Z^2\frac{\varepsilon\tilde{\varepsilon}} {\varepsilon-\tilde{\varepsilon}}}))^{C_{p,\tilde{\varepsilon}}}
\end{align*}
for some function $C_{p,\cdot}:(p^2\kappa,\varepsilon)\to[0,\infty)$.
\end{thm}
\begin{proof} By \eqref{em} for $\tilde{\varepsilon}\in(\kappa,\varepsilon)$ and \cite[Theorem 5.2]{Wp}, we conclude that $\mu^Z$ is the unique invariant probability measure of \eqref{EMI} in $\scr P_Z$ and $\rho$ has a strictly positive continuous version. So, it remains to prove the three estimates above.
By \cite[Theorem 5.2]{Wp}, we derive $\log \rho\in H^{1,2}_{\sigma}(\mu^0)$ with
\begin{align*}
\mu^0(|\sigma^\ast\nabla\log \rho|^2)\leq 4\mu^0(|Z(\cdot,\mu^Z)|^2)\leq 4\mu^0((|Z(\cdot,\mu^0)|+2K_Z)^2),
\end{align*}
and for any $\tilde{\varepsilon}\in(\kappa,\varepsilon)$,
\begin{align*}
\nonumber &\mu^0(|\sigma^\ast\nabla\sqrt{\rho}|^2)\leq \frac{1}{\tilde{\varepsilon}-\kappa}(\log \mu^0(\e^{\tilde{\varepsilon } |Z(\cdot,\mu^Z)|^2})+\beta)\leq \frac{1}{\tilde{\varepsilon}-\kappa}(\log \mu^0(\e^{ \varepsilon|Z(\cdot,\mu^0)|^2+4 K_Z^2\frac{\varepsilon\tilde{\varepsilon}} {\varepsilon-\tilde{\varepsilon}}})+\beta),
\end{align*}
here we used \eqref{em} in the last step.

Next, for any $p\in(1,\frac{\sqrt{\varepsilon}}{\sqrt{\kappa}})$ and $\tilde{\varepsilon}\in(p^2\kappa,\varepsilon)$, we have $p\in(1,\frac{\sqrt{\tilde{\varepsilon}}}{\sqrt{\kappa}})$. Again by \eqref{em} for $\tilde{\varepsilon}\in(p^2\kappa,\varepsilon)$ and \cite[Theorem 5.2]{Wp}, there exists a  constant $C_{\tilde{\varepsilon},p}\geq 1$ such that
\begin{align*}
\mu^0(|\sigma^\ast\nabla\rho^{\frac{p}{2}}|^2+\rho^p)\leq C_{p,\tilde{\varepsilon}}(\mu^0(\e^{\tilde{\varepsilon}|Z(\cdot,\mu^Z)|^2}))^{C_{p,\tilde{\varepsilon}}}\leq C_{p,\tilde{\varepsilon}}(\mu^0(\e^{ \varepsilon|Z(\cdot,\mu^0)|^2+4 K_Z^2\frac{\varepsilon\tilde{\varepsilon}} {\varepsilon-\tilde{\varepsilon}}}))^{C_{p,\tilde{\varepsilon}}}.
\end{align*}
So, the proof is completed.
\end{proof}
\section{Existence and Uniqueness of Invariant Probability Measure}
In this section, we will consider two cases: one is the symmetric case and the other one is stochastic Hamiltonian system. We will consider the invariant probability measure of \eqref{E1} in the class:
$$\tilde{\scr P}_Z=\{\nu=\rho_\nu\mu^0:\nu \ \ \text{is an invariant probability measure of \eqref{E1}}\}.$$
For any $\mu\in\scr P(\R^d)$, consider
\beq\label{EMO} \d X_t=\{ Z_0(X_t)+\sigma(X_t)Z(X_t, \mu)\}\d t +\si(X_t)\d W_t,\end{equation}
and denote
$$\scr P^\mu_Z=\{\nu=\rho_\nu\mu^0:\nu \ \ \text{is an invariant probability measure of \eqref{EMO}}\}.$$
\subsection{Symmetric Case}
\begin{exa}\label{non}  Let $Z(x,\mu)=\frac{\sqrt{2}}{2}(\nabla F(x,\mu)+\nabla \bar{F}(x))$ and $\sigma=\sqrt{2}I_{d\times d}$. Assume    $({\bf H})$ and that there exist constants $\varepsilon>\kappa, C>0$ and $\delta\in(0,\frac{\log{2}}{2})$ such that
$$\mu^0(\e^{F(\cdot,\mu^0)+\bar{F}})+\mu^0(\e^{\frac{\varepsilon}{2}|\nabla F(\cdot,\mu^0)+\nabla \bar{F}|^2})<\infty,$$
$$|\nabla F(x,\mu)-\nabla F(x,\nu)|\leq C\|\mu-\nu\|_{TV},\ \ \mu,\nu\in\scr P(\R^d), x\in\R^d,$$
and $$| F(x,\mu)-F(x,\nu)|\leq\delta\|\mu-\nu\|_{TV},\ \ \mu,\nu\in\scr P(\R^d), x\in\R^d.$$
Then \eqref{E1} has a unique invariant probability measure in $\tilde{\scr P}_Z$.
\end{exa}
\begin{proof} By \eqref{em} for $\mu^Z=\mu, Z(x,\mu)=\frac{\sqrt{2}}{2}(\nabla F(x,\mu)+\nabla \bar{F}(x))$ and $\tilde{\varepsilon}\in(\kappa,\varepsilon)$, \cite[Theorem 5.2]{Wp} implies that for any $\mu\in\scr P(\R^d)$, \eqref{EMO} has a unique invariant probability measure in $\scr P^\mu_Z$, which is denoted by $\Gamma(\mu)$.
Therefore, $\Gamma$ construct a map from $\scr  P(\R^d)$ to $\scr  P(\R^d)$. Moreover,
it is clear that
$$\frac{\d\Gamma(\mu)}{\d \mu^0}=\frac{\e^{F(\cdot,\mu)+\bar{F}}}{\mu^0(\e^{F(\cdot,\mu)+\bar{F}})}.$$
By Taylor's expansion, we arrive at
\begin{align*} \left|\e^{F(x,\mu)}-\e^{F(x,\nu)}\right|&\leq \e^{F(x,\nu)}\sum_{k=1}^\infty\frac{|F(x,\mu)-F(x,\nu)|^k}{k!}\\
&\leq \e^{F(x,\nu)}\sum_{k=1}^\infty\frac{\delta^k\|\mu-\nu\|_{TV}^k}{k!}\\
&\leq \e^{F(x,\nu)}\|\mu-\nu\|_{TV}\sum_{k=1}^\infty\frac{\delta^k2^{k-1}}{k!}\\
&=\e^{F(x,\nu)}\|\mu-\nu\|_{TV}\frac{\e^{2\delta}-1}{2}.
\end{align*}
As a result, it holds
\begin{align*}
&\|\Gamma(\mu)-\Gamma(\nu)\|_{TV}\\
&=\int_{\R^d}\left|\frac{\e^{F(\cdot,\mu)+\bar{F}}}{\mu^0(\e^{F(\cdot,\mu)+\bar{F}})} -\frac{\e^{F(\cdot,\nu)+\bar{F}}}{\mu^0(\e^{F(\cdot,\nu)+\bar{F}})}\right|\mu^0(\d x)\\
&= \int_{\R^d}\left|\frac{\e^{F(\cdot,\mu)+\bar{F}}\mu^0(\e^{F(\cdot,\nu)+\bar{F}})-\e^{F(\cdot,\nu)+\bar{F}}\mu^0(\e^{F(\cdot,\mu)+\bar{F}})} {\mu^0(\e^{F(\cdot,\mu)+\bar{F}})\mu^0(\e^{F(\cdot,\nu)+\bar{F}})}\right|\mu^0(\d x)\\
&\leq \int_{\R^d}\frac{\left|\mu^0(\e^{F(\cdot,\nu)+\bar{F}})-\mu^0(\e^{F(\cdot,\mu)+\bar{F}})\right|}{\mu^0(\e^{F(\cdot,\nu)+\bar{F}})}\mu^0(\d x)+\int_{\R^d}\frac{\left|\e^{F(\cdot,\mu)+\bar{F}}-\e^{F(\cdot,\nu)+\bar{F}}\right|}{\mu^0(\e^{F(\cdot,\nu)+\bar{F}})}\mu^0(\d x)\\
&\leq \|\mu-\nu\|_{TV}(\e^{2\delta}-1)
\end{align*}
So, when $\delta\in(0,\frac{\log{2}}{2})$, $\Gamma$ is a strictly contractive map on $(\scr P(\R^d),\|\cdot\|_{TV})$. By Banach's fixed point theorem, we finish the proof.
\end{proof}
\subsection{Stochastic Hamiltonian system}
The next example concentrates on the stochastic Hamiltonian system in $\R^{2d}$.
\begin{exa}\label{SHS} Let $Z_0(x,y)=(y,-x-y), x,y\in\R^d$, $\sigma=\left(
                                                                      \begin{array}{cc}
                                                                        0_{d\times d} & 0_{d\times d} \\
                                                                        0_{d\times d}  &\sqrt{2}I_{d\times d}  \\
                                                                      \end{array}
                                                                    \right)
$. Consider
\begin{align}\label{SHS00}\left\{
  \begin{array}{ll}
    \d X_t=Y_t \d t\\
    \d Y_t=(-X_t-Y_t)\d t+\nabla H(\cdot,\L_{(X_t,Y_t)})(X_t)\d t+\nabla \bar{H}(X_t)\d t+\sqrt{2}\d W_t,
  \end{array}
\right.
\end{align}
here $H:\R^d\times \scr P(\R^{2d})\to\R$, $\bar{H}:\R^d\to\R$. Let $Z(x,\mu)=\frac{\sqrt{2}}{2}(0, \nabla H(\cdot,\mu)(x)+\nabla \bar{H}(x))$. By \cite[Example 5.1]{W17}, $\mu^0(\d x,\d y)=\frac{1}{(2\pi)^d}\e^{-\frac{|x|^2+|y|^2}{2}}\d x\d y$ and $$\|P_{t_0}^0\|_{L^2(\mu^0)\to L^{4}(\mu^0)}=1$$ for some $t_0>0$. Assume that there exist constants $\varepsilon>\frac{5}{2}t_0, C>0$ and $\delta\in(0,\frac{\log{2}}{2})$ such that
$$\mu^0(\e^{H(\cdot,\mu^0)+\bar{H}})+\mu^0(\e^{\frac{\varepsilon}{2}|\nabla H(\cdot,\mu^0)+\nabla \bar{H}|^2})<\infty,$$
$$|\nabla H(x,\mu)-\nabla H(x,\nu)|\leq C\|\mu-\nu\|_{TV},\ \ \mu,\nu\in\scr P(\R^{2d}), x\in\R^d,$$
and $$| H(x,\mu)-H(x,\nu)|\leq\delta\|\mu-\nu\|_{TV},\ \ \mu,\nu\in\scr P(\R^{2d}), x\in\R^d.$$
Then \eqref{SHS00} has a unique invariant probability measure in $\tilde{\scr P}_{Z}$.
\end{exa}
\begin{proof} Note that \eqref{HYP} holds for $p_0=2$ and so $\kappa_0=\frac{5}{2}t_0$. Using \eqref{em} for $\mu^Z=\mu, Z(x,\mu)=\nabla H(x,\mu)+\nabla \bar{H}(x)$ and $\tilde{\varepsilon}\in(\kappa_0,\varepsilon)$ and applying \cite[Theorem 3.1, Theorem 4.1]{Wp}, we conclude that \eqref{EMO} has a unique invariant probability measure in $\scr P^\mu_Z$ and we denote it as $\Gamma(\mu)$, which satisfies
$$\frac{\d\Gamma(\mu)}{\d \mu^0}=\frac{\e^{H(\cdot,\mu)+\bar{H}}}{\mu^0(\e^{H(\cdot,\mu)+\bar{H}})}.$$
In fact, the infinitesimal generator of \eqref{SHS00} with $\mu$ in place of $\L_{(X_t,Y_t)}$ is
$$L f(x,y)=y\nabla _x f+(-x-y+\nabla H(x,\mu)+\nabla \bar{H}(x))\nabla_y f+\nabla^2 _y f, \ \ f\in C_0^\infty(\R^{2d}),$$
which yields
\begin{align*}&\int_{\R^{2d}} L f(x,y)\exp\left\{-\frac{|x|^2}{2}-\frac{|y|^2}{2}+H(x,\mu)+\bar{H}(x)\right\}\d x\d y\\
&=-\int_{\R^{2d}}\left\<\nabla _xf, \nabla_y\exp\left\{-\frac{|x|^2}{2}-\frac{|y|^2}{2}+H(x,\mu)+\bar{H}(x)\right\}\right\>\d x\d y\\
&+\int_{\R^{2d}}\left\<\nabla _yf, \nabla_x\exp\left\{-\frac{|x|^2}{2}-\frac{|y|^2}{2}+H(x,\mu)+\bar{H}(x)\right\}\right\>\d x\d y=0,\ \ f\in C_0^\infty(\R^{2d}).
\end{align*}
Then repeating the proof of Example \ref{non}, we complete the proof.
\end{proof}
  \paragraph{Acknowledgement.} The authors   would like to thank Professor Feng-Yu Wang for helpful comments.

\end{document}